\documentclass[]{article}
\input{amssym}
\input{epsf}
\begin{document}
\def\be{\begin{eqnarray*}}
\def\ee{\end{eqnarray*}}
\def\di{\displaystyle}
\title{{\bf Symmetry Analysis of Cylindrical Helmholtz Equation}}
\author{ Mehdi Nadjafikhah\thanks{Department of
Mathematics, Iran University of Science and Technology, Narmak-16,
Tehran, I.R.Iran. e-mail:~m\_nadjafikhah@iust.ac.ir}\and Ali
Mahdipour--Shirayeh\thanks{e-mail:~mahdipour@iust.ac.ir}}
\date{ }
\maketitle
\begin{abstract}
In this paper, we present the point symmetry group of
three-dimensional homogeneous Helmholtz equation, when we consider
the cylindrical coordinate system. In continuation, we present a
complete set of functionally independent invariants of the
equation along with the form of the general solution provided by
these invariants. Finally, we find an optimal system of
one-dimensional Lie subgroups of the full symmetry group.
\end{abstract}

\medskip \noindent {\bf A.M.S. 2000 Subject Classification:} 34C14, 35J05, 70G65.


%
%*****************************************
\section{Introduction}
The Helmholtz equation, named for Hermann von Helmholtz, is the
elliptic partial differential equation (in homogenous case)
\begin{eqnarray*}
\nabla\,u+k^2u=0,
\end{eqnarray*}
when $u$ is a vector function and $\nabla$ is the vector
Laplacian.\par
The three-dimensional homogeneous {\it cylindrical Helmholtz
equation} (in brief CHE), is a Helmholtz differential equation
defined in the cylindrical $(r,\theta,z)$ coordinate system, i.e.
\begin{eqnarray}
\textrm{CHE}\;:\;u_{rr}+\frac{1}{r}u_r+\frac{1}{r^2}u_{\theta\theta}+u_{zz}+k^2u=0,\label{eq:1}
\end{eqnarray}
where $u$ is a function of variables $r,\theta,z$ and $k$ is a
constant.\par
The Helmholtz equation naturally appears from general conservation
laws of physics and often arises in the study of physical problems
involving partial differential equations (PDEs) in both space and
time, more particularly in mechanics, acoustics, electromagnetics,
etc. (see for details \cite{FK,GD,Ho}). \par
This equation can be interpreted as a wave equation for
monochromatic waves and can be derived from the heat conduction
equation, Schr{\"{o}}dinger equation, telegraph and other wave
type, or evolutionary equations \cite{GD}. When $k=0$, the
Helmholtz differential equation reduces to Laplace's equation and
when $k^2<0$ (for imaginary), the equation becomes the space part
of the diffusion equation.\par
In the mathematical sense, it appears as an eigenvalue problem for
the Laplace operator $\nabla^2$ (see \cite{FK,GD,MS}). As is
well-known, there are some mathematical theory which can be
applied to a solution of the Helmholtz equation in Cartesian,
spherical, cylindrical coordinates, and etc., by applying the
technique of separation of variables and boundary element method,
its various solutions for many basic shapes, for example in a
vibrating membrane, along with other different methods
\cite{AK,BG,GD,MS}. In addition, since the Helmholtz equation is
an elliptic equation, so it is usual to consider boundary value
problem along with transmission conditions and integral
representation of solutions (see \cite{FK,GD} and references
therein).\par
At the present work, we are dealing with the symmetry analysis and
a one-dimensional optimal system of the Helmholtz equation in
cylindrical coordinate. Cylindrical coordinates are the natural
choice for studying different shapes of the Helmholtz equation,
that is a good motivation for work with this coordinate. Also we
find a form of general solutions providing by the differential
invariants of the equation.
%
%********************************************
\section{Symmetry of CHE}
In this section, we are deal with symmetry analysis of CHE, when
we act the point symmetry group on the second order jet space
$J^2=J^2({\Bbb R}^3,{\Bbb R})$ with local coordinates
$(r,\theta,z,u,u^{(1)},u^{(2)})$. In these charts $u^{(\alpha)}$
($\alpha=1,2$) indicates the set of all $\alpha^{th}$ derivatives
of the dependent variable $u$ with respect to independent variable
$r,\theta,z$.\par
The base space of the jet bundle $J^2$ is the space ${\cal M} =
X\times U$ of independent and dependent variables $(x,u)$ for
$x=(r,\theta,z)$, that is in fact the zero-order jet space
$J^0({\Bbb R}^3,{\Bbb R})$. CHE defines a subvariety of jet bundle
$J^2$ of three-dimensional submanifold of ${\cal M}$, that is,
graph of function $u=u(r,\theta,z)$.\par
The point symmetry group on ${\cal M}$ is defined by
transformations in the general form of
\begin{eqnarray*}
\widetilde{r}=\varphi(r,\theta,z,u),\hspace{0.7cm}\widetilde{\theta}=\psi(r,\theta,z,u),\hspace{0.7cm}
\widetilde{z}=\chi(r,\theta,z,u),\hspace{0.7cm}\widetilde{u}=\zeta(r,\theta,z,u),
\end{eqnarray*}
for arbitrary smooth functions $\varphi,\psi,\chi$ and
$\zeta$.\par
Let
\begin{eqnarray}
v=\sum_{i=1}^3\,\xi^i(x,u)\frac{\partial}{\partial
x_i}+\eta(x,u)\frac{\partial}{\partial u},\label{eq:2}
\end{eqnarray}
be the general form of an infinitesimal generator of point
symmetry transformations on ${\cal M}$, in which the coefficients
are arbitrary smooth functions. Let the second order prolongation
of $v$ be defined as
\begin{eqnarray*}
v^{(2)}= v +
\sum_{\#J=j=1}^2\,\widehat{\eta}_J(x,u^{(j)})\frac{\partial}{\partial
u_J},
\end{eqnarray*}
where we assumed that each index signifies to the derivation in
respect to the index, and $J=(j_1,j_2)$ is a multi-index with
length $0\leq\#J\leq 2$. In addition, coefficients are given by
the prolongation formula (see \cite{Ol})
\begin{eqnarray*}
\widehat{\eta}_J=D_J(Q)+\sum_{i=1}^3\xi^i\,u_{J,i}.
\end{eqnarray*}
In the latter relation, $D_J$ is the total derivative with respect
to the variables signified by $J$, and $Q=\eta -
\sum_{i=1}^3\xi^i\,u_i$ is the characteristic of $v$.\par
%More explicitly, for $1\leq i,j\leq 3$ we have
%
%\begin{eqnarray*}
%\hspace{-0.7cm} \widehat{\eta}_i &=& \eta_i + \eta_u\,u_i-\sum_{k=1}^3(\xi^k_i u_i + \xi^k_u u_i u_k) ,\\
%\hspace{-0.7cm} \widehat{\eta}_{ij} &=& \eta_{ij} + \eta_{i u}u_j
%+ \eta_{u j}u_i + \eta_{u u}u_i u_j - \sum_{k=1}^3(\xi^k_{ij}u_k +
%\xi^k_{i u}u_j u_k \\
%&&+ \xi^k_{u j}u_i u_k  + \xi^k_{uu}u_i u_j u_k+ \xi^k_u u_i
%u_{kj} + \xi^k_u u_k u_{ij} + \xi^k_{j}u_{ki} + \xi^k_u u_j
%u_{ki}).
%\end{eqnarray*}
%
According to \cite{Ol}, the vector field $v$ is an infinitesimal
symmetry of Eq. (\ref{eq:1}), if and only if it satisfies the
infinitesimal invariance condition
\begin{eqnarray}
v^{(2)}\Big(u_{rr}+\frac{1}{r}u_r+\frac{1}{r^2}u_{\theta\theta}+u_{zz}+k^2u\Big)=0,\label{eq:3}
\end{eqnarray}
whenever Eq. (\ref{eq:1}) vanishes.\par
Eq. (\ref{eq:3}) leads to the
infinitesimal determining equation, that must be solved for
arbitrary functions $\xi^i$ and $\eta$, when they are dependent to
variables $r,\theta,z,u$ and not to the derivatives $u$ in respect
to its variables. Hence, the induced equation satisfies for
$\xi^i$ and $\eta$, if and only if the following determining
equations are fulfilled
\begin{eqnarray}
&&\hspace{-1cm} \xi^1_{uu}=0,\hspace{1.8cm}  \xi^2_{uu}=0,\hspace{1.5cm}  \xi^2_{uu}=0, \hspace{1.3cm}\xi^2_u=0, \nonumber\\[1mm]
&&\hspace{-1cm} \xi^1_{uu}=0,\hspace{1.8cm}
\xi^3_{uu}=0,\hspace{1.5cm}
\xi^1_u=0, \hspace{1.5cm}  \xi^1_u=0, \nonumber\\[1mm]
&&\hspace{-1cm} \xi^2_u=0,\hspace{1.97cm}
\xi^1_u=0,\hspace{1.7cm}  \xi^3_u=0,\hspace{1.5cm} \xi^3_u=0, \nonumber\\[1mm]
&&\hspace{-1cm} r^2\xi^2_z + \xi^3_{\theta}=0,\hspace{0.9cm}
\xi^1_u=0,\hspace{1.7cm} \xi^3_{uu}=0, \hspace{1.35cm}\xi^1_z + \xi^3_r=0, \nonumber\\[1mm]
&&\hspace{-1cm} 2\,\xi^1_u +
r\eta_{uu}-2\,r\xi^1_{ru}=0,\hspace{0.85cm}
r\xi^1_r-\,r\xi^2_{\theta}-\xi^1=0, \hspace{0.9cm} \xi^1_r-\xi^3_z=0, \nonumber\\[1mm]
&&\hspace{-1cm} 2\,\xi^3_{zu}-\eta_{u u}=0,\hspace{2.17cm}
r^2\xi^2_r+\xi^1_{\theta}=0,\hspace{1.7cm} \xi^3_{ru}+\xi^1_{zu}=0,\label{eq:5}\\[1mm]
&&\hspace{-1cm} \xi^3_{\theta u}+r^2\xi^2_{zu}=0,\hspace{2.1cm}
\eta_{uu}-2\,\xi^2_{\theta u}=0, \hspace{1.5cm} r^2\xi^2_{ru}+\xi^1_{\theta u}=0, \nonumber\\[1mm]
&&\hspace{-1cm} r^2\,k^2 u\,\xi^2_u + 2\,\eta_{\theta u}-r\xi^2_r
-\xi^2_{\theta\theta}-r^2\xi^2_{rr}-r^2\xi^2_{zz}=0,  \nonumber\\[1mm]
&&\hspace{-1cm}  3\,r^2\,k^2u\,\xi^1_u +
r\,\xi^1_r+2\,r^2\eta_{ru}-r^2\xi^1_{zz}
-\xi^1 -\xi^1_{\theta\theta}-r^2\xi^1_{rr}=0, \nonumber\\[1mm]
&&\hspace{-1cm}
r^2\,k^2u\,\xi^3_u+2\,r^2\eta_{zu}-r\xi^3_r-r^2\xi^3_{zz}
-\xi^3_{\theta\theta}-r^2\xi^3_{rr}=0, \nonumber\\[1mm]
&&\hspace{-1cm}  r^2\eta_{z,z}+2\,r^2\,k^2 u(\xi^1_r - \eta_u)
+r^2\eta_{rr}+r\eta_r + \eta_{\theta\theta} + r^2\eta\,k^2 =0.
\nonumber
\end{eqnarray}
The solution space ${\goth g}$ of the system of determining
equations (\ref{eq:5}) with respect to the unknown coefficients of
$v$, is the Lie algebra of infinitesimal symmetry of Eq.
(\ref{eq:1}). By solving system (\ref{eq:5}) (that is applicable
with \textsc{Maple}) we obtain the below solutions
\begin{eqnarray}
\xi^1 &=& (c_1\,z+ c_3)\sin\theta+(c_2\,z +c_4)\cos\theta, \nonumber\\
\xi^2 &=& -\frac{1}{r}(c_2\,z+c_4)\sin\theta+\frac{1}{r}(c_1\,z+c_3)\cos\theta+ c_5,\nonumber\\
\xi^3 &=& -r(c_1\,\sin\theta + c_2\,\cos\theta) + c_6, \label{eq:6}\\
\eta &=&
\exp(-\sqrt{c_3}z-\sqrt{c_2}\theta)(c_7\exp(2\sqrt{c_2}\theta)
+c_8)(c_9\exp(2\sqrt{c_3}z)+c_{10})  \nonumber\\
&&\times\left(c_{11}\,{\rm BY}
(\sqrt{-c_2},\sqrt{c_3+k^2}r)+c_{12}\,{\rm
BJ}(\sqrt{-c_2},\sqrt{c_3+k^2}r)\right) + c_{13}\,u, \nonumber
\end{eqnarray}
where, ${\rm BJ}$ and ${\rm BY}$ are the Bessel functions of the
first and second kinds, resp. which satisfy Bessel's equation
$x^2\,y''+x\,y'+(x^2-v^2)y=0$, and $c_i$ (for $1\leq i\leq 13$)
are arbitrary constants.\par
Noting to the concluded form of infinitesimal operators of Eq.
(\ref{eq:1}), one may divide the resulting form of $v$ to the set
of following vector fields, as a basis of Lie algebra ${\goth g}$
of the symmetry group $G$ of CHE
\begin{eqnarray}
&&\hspace{-1.5cm} X_1 =
\frac{\partial}{\partial\theta},\hspace{1.5cm} X_2 =
\frac{\partial}{\partial z},\hspace{1.4cm}
X_3 = u\,\frac{\partial}{\partial u}, \nonumber\\
&&\hspace{-1.5cm} X_4 = \sin\theta\,\frac{\partial}{\partial
r}+\frac{\cos\theta}{r}\,\frac{\partial}{\partial\theta},\hspace{2cm}
X_5 =  -\cos\theta\,\frac{\partial}{\partial
r}+\frac{\sin\theta}{r}\,\frac{\partial}{\partial\theta},\nonumber\\[-2mm]
&&\label{eq:7}\\[-2mm]
&&\hspace{-1.5cm} X_6 = -z\,\cos\theta\,\frac{\partial}{\partial
r}+\frac{z\,\sin\theta}{r}\,\frac{\partial}{\partial\theta}
+r\,\cos\theta\,\frac{\partial}{\partial z},\nonumber\\
&&\hspace{-1.5cm} X_7 = z\,\sin\theta\,\frac{\partial}{\partial
r}+\frac{z\,\cos\theta}{r}\,\frac{\partial}{\partial\theta}
-r\,\sin\theta\,\frac{\partial}{\partial z}.\nonumber
\end{eqnarray}
The commutator of every $X_i$ and $X_j$ i.e. $[X_i,X_j]:=X_i\,X_j
- X_j\,X_i$ ($1\leq i,j\leq7$), is straightforwardly a linear
combination of vector fields (\ref{eq:7}). The complete set of
commutators is given in Table \ref{tabel:1}.\\
\begin{table}[h]
\begin{center}
\begin{tabular}{|c|ccccccc|}
\hline  &       &       &      &       &       &       &       \\[-2mm]
        & $X_1$ & $X_2$ & $X_3$& $X_4$ & $X_5$ & $X_6$ & $X_7$ \\[1mm]
\hline  &       &       &      &       &       &       &       \\[-2mm]
  $X_1$ & 0     & 0     & 0    & $-X_5$& $X_4$ & $X_7$ &$-X_6$ \\[1mm]
  $X_2$ & 0     & 0     & 0    & 0     & 0     & $X_5$ & $X_4$ \\[1mm]
  $X_3$ & 0     & 0     & 0    & 0     & 0     & 0     & 0     \\[1mm]
  $X_4$ & $X_5$ & 0     & 0    & 0     & 0     & 0     &$-X_2$ \\[1mm]
  $X_5$ & $-X_4$& 0     & 0    & 0     & 0     & $-X_2$& 0     \\[1mm]
  $X_6$ & $-X_7$& $-X_5$& 0    & 0     & $X_2$ & 0     & $X_1$ \\[1mm]
  $X_7$ & $X_6$ & $-X_4$& 0    & $X_2$ & 0     & $-X_1$& 0     \\[1mm]
\hline
\end{tabular}
\caption{Commutator table for ${\goth g}$.} \label{tabel:1}
\end{center}
\end{table}
\paragraph{Theorem 1.}{\it The complete set of infinitesimal generators of
point symmetry group of Eq. (\ref{eq:1}), is introduced in
(\ref{eq:2}) and (\ref{eq:6}) and is a seven-dimensional Lie
algebra with commutators given in Table \ref{tabel:1}.}\\

Since the derived algebra ${\goth g}^{(1)}=[{\goth g},{\goth
g}]={\goth g}\diagup\langle X_3 \rangle$ and for each $k>1$,
${\goth g}^{(k+1)}=[{\goth g}^{(k)},{\goth g}^{(k)}]={\goth
g}^{(1)}$, so the descending sequence of the derived subalgebras
of ${\goth g}$,
\begin{eqnarray*}
{\goth g}\supset{\goth g}^{(1)}={\goth g}^{(2)}=\cdots
\end{eqnarray*}
does not terminate with a null ideal and hence ${\goth g}$ is not
solvable \cite{Ly,Ov}. Also, ${\goth g}$ is not semisimple, since
by applying Cartan's criterion of semisimplicity \cite{Ov}, its
Killing form defined for $v=\sum_i v_i X_i$ and $w=\sum_j w_j X_j$
in ${\goth g}$ as
\begin{eqnarray*}
K(v,w)={\rm tr}({\rm ad}(v)\circ{\rm ad}(w))=-4(v_1 w_1 + v_6 w_6
+ v_7w_7)
\end{eqnarray*}
is degenerate.\par
The radical of ${\goth g}$ is equal to the Lie subalgebra ${\goth
r}:=\langle X_2,X_3,X_4,X_5\rangle$ and hence Levi decomposition
of ${\goth g}$ is in the form of ${\goth g}={\goth r}\oplus{\goth
s},$ for semisimple subalgebra ${\goth s}:={\goth g}\diagup{\goth
r}=\langle X_1,X_6,X_7 \rangle$. According to the Table
\ref{tabel:1}, ${\goth r}$ is an abelian group of dimension 4,
thus is isomorphic to ${\Bbb R}^4$, and the (real) Lie subalgebra
${\goth s}$ is isomorphic to the Lie algebra ${\goth s}{\goth
o}(3)$ of the special orthogonal group. In addition, ${\goth s}$
is not an ideal of ${\goth g}$ and therefore we find the structure
of the symmetry algebra of CHE as following semi-direct product
\begin{eqnarray}
{\goth g}={\Bbb R}^4\ltimes{\goth s}{\goth o}(3).\label{eq:8}
\end{eqnarray}
For finding the symmetry group $G$ of the Lie algebra ${\goth g}$
generated by infinitesimal operators (\ref{eq:7}), we determine
the one-parameter group (flow) of these vector fields:
\paragraph{Theorem 2.}{\em If we denote the group transformations of
each infinitesimal generator $X_i$ ($i=1,\cdots,7$) with $g_i(s)$
for different values of parameter $s$, then we have
\begin{eqnarray*}
g_1(s)&:&(r,\theta,z,u)\longmapsto(r,\theta + s ,z,u),\\
g_2(s)&:&(r,\theta,z,u)\longmapsto(r,\theta,z + s ,u),\\
g_3(s)&:&(r,\theta,z,u)\longmapsto(r,\theta ,z,u + s),\\
g_4(s)&:&(r,\theta,z,u)\longmapsto \Bigg(r\Big(\cos^2\theta +
\frac{1}{4r^2}\left(2s-r\sin 2\theta\right)^2\Big)^{1/2},\\
&&\arctan\left(\frac{s}{r}(1+\tan^2\theta)^{1/2}-\tan \theta\right),z,u\Bigg),\\
g_5(s)&:&(r,\theta,z,u)\longmapsto
\Bigg(\Big(s-\!\frac{r}{(1\!+\!\tan^2\theta)^{1/2}}\Big)\!
\Big(1\!+\!\Big(\!\frac{r\tan\theta}{s(1\!+\!\tan^2\theta)^{1/2}\!-\!r}\Big)^2\Big)^{1/2}\!\!,\\
&&\arctan\Big(\!\frac{r\,\tan\theta}{s(1\!+\!\tan^2\theta)^{1/2}\!+\! r}\Big)\!,z,u\Bigg),\\
g_6(s)&:&(r,\theta,z,u)\longmapsto \Bigg(\!\Big(\!c_1\cos^2(s+c_2)\!+\!d_1^2\sin^2(s+c_2)\Big)^{1/2}\!\!\!\!,\\
&&\theta\!-\!\frac{\Big((c_1-d_1^2)\cos^2(s+c_2)\Big)^{1/2}\!\!\!\!\tanh^{-\!1}\!\!\Big(\frac{d_1^2-c_1}{d_1}\!\cos(s+c_2)\!\Big)
}{(d_1^2-c_1)^{1/2}\cos(s+c_2)}\\
&&+\frac{\Big((c_1-d_1^2)\cos^2(c_2)\Big)^{1/2}\!\!\!\!\tanh^{-1}\Big(\frac{(d_1^2-c_1)}{d_1}\cos(c_2)\Big)}
{(d_1^2-c_1)^{1/2}\cos(c_2)},\\
&& \frac{\sqrt{2}}{2}\Big((c_1-d_1^2)\sin(2\,s+2\,c_2)\Big)^{1/2},u \Bigg),\\
g_7(s)&:&(r,\theta,z,u)\longmapsto \Bigg(\!\Big(\!c_1\cos^2(s+c_3)\!+\!d_2^2\sin^2(s+c_3)\Big)^{1/2}\!\!\!\!,\\
&&\theta\!+\!\frac{\Big((c_1-d_2^2)\cos^2(s+c_3)\Big)^{1/2}\!\!\!\!\tanh^{-\!1}\!\!\Big(\frac{d_2^2-c_1}{d_1}\!\cos(s+c_3)\!\Big)
}{(d_2^2-c_1)^{1/2}\cos(s+c_3)}\\
&&-\frac{\Big((c_1-d_2^2)\cos^2(c_3)\Big)^{1/2}\!\!\!\!\tanh^{-1}\Big(\frac{(d_2^2-c_1)}{d_1}\cos(c_3)\Big)}
{(d_2^2-c_1)^{1/2}\cos(c_3)},\\
&&\frac{\sqrt{2}}{2}\Big((c_1-d_2^2)\sin(2\,s+2\,c_3)\Big)^{1/2},u
\Bigg),
\end{eqnarray*}
where we assume that $c_1=r^2 + z^2$, $c_2 = \arctan( z(
r^2\,\cos^2 t)^{-1/2})$, $c_3=-s-\arctan(z(r^2\,\sin^2t)^{-1/2})$,
$d_1=r\,\sin t$ and , $d_2=r\,\cos t$.}\\

It is well-known that for each one-parameter subgroup of the full
symmetry group of a system there will correspond a family of
solutions that are called invariant solutions \cite{Ol,OL,Ov}. In
section four, we introduce a general solution of CHE with respect
to its functionally independent invariants.
\paragraph{Theorem 3.}
{\em If $u=f(r,\theta,z)$ be a solution of (\ref{eq:1}), then are
the below functions ($i=4,\cdots,7$)\vspace{-0.4cm}
\begin{eqnarray*}
&&g_1(s)\cdot f(r,\theta,z)=f(r,\theta + s ,z),\nonumber\\
&&g_2(s)\cdot f(r,\theta,z)=f(r,\theta,z + s ),\nonumber\\
&&g_3(s)\cdot f(r,\theta,z)=f(r,\theta ,z)-s,\nonumber\\
&&g_i(s)\cdot f(r,\theta,z)=f(g_i(s)\cdot(r,\theta,z)),\nonumber
\end{eqnarray*}
when $g_i(s)\cdot(r,\theta,z)$ means the first three components of
$g_i(s)\cdot(r,\theta,z,u)$ and $s$ is an arbitrary parameter.}
%
%******************************************
\section{Invariants and General Solutions}
In this section, we study the invariant solutions of various
combinations of vector fields (\ref{eq:7}) based on the following
method:

the invariants $u=I(r,\theta,z,u)$ of one--parameter group with
infinitesimal generators in the form of (\ref{eq:7}) satisfy the
linear homogeneous partial differential equations of first order:
\begin{eqnarray*}
v[I]=0.
\end{eqnarray*}
The solutions of the latter are found by the method of
characteristics (See \cite{Ib} and \cite{OL} for details). So we
can replace the latest Eq. by the following characteristic system
of ordinary differential equations
\begin{eqnarray}
\frac{dr}{\xi_1}=\frac{d\theta}{\xi_2}=\frac{dz}{\xi_3}=\frac{du}{\eta}
.\label{eq:9}
\end{eqnarray}
By solving the Eqs. (\ref{eq:9}) of differential generators
(\ref{eq:6}), we (locally) find the following general solutions
\begin{eqnarray}
I_1(r,\theta,z,u) &=& 2z +\frac{(c_1\,\sin\theta +
c_2\,\cos\theta)\,r^2 + c_6\,r}{(c_1\,z+ c_3)\sin\theta+(c_2\,z
+c_4)\cos\theta}=d_1 , \nonumber\\
I_2(r,\theta,z,u) &=& z - \frac{1}{c_{13}}(c_1\,\sin\theta +
c_2\,\cos\theta)\,\ln\Big((c_{11}\,{\rm
BY}(\sqrt{-c_2},\sqrt{c_3+k^2})  \nonumber\\
&&+ c_{12}\,{\rm
BJ}(\sqrt{-c_2},\sqrt{c_3+k^2}))(c_7c_9\exp(\sqrt{c_2}\,\theta+\sqrt{c_3}\,z)
\nonumber\\[2mm]
&& + c_7c_{10}\exp(\sqrt{c_2}\,\theta-\sqrt{c_3}\,z)
+ c_8 c_9\exp(-\sqrt{c_2}\,\theta+\sqrt{c_3}\,z)  \nonumber\\[-2mm]
&& \label{eq:10}\\[-2mm]
&& + c_8 c_{10}\exp(-\sqrt{c_2}\,\theta-\sqrt{c_3}\,z) )+c_{13}\,u
\Big)= d_2 , \nonumber\\
I_3(r,\theta,z,u) &=&
-\frac{1}{B}\arctan(\alpha)\left(2\,c_5\,A^{-1}\,(c_2c_3-c_1c_4)r^3
+ 2c_6\,r\right) \nonumber\\
&&+r^2\,A^{-1}\,\Big(\ln(C)-\ln\Big(\frac{1}{\cos\theta+1}\Big)\Big)(c_2^2\,z
+ c_1c_3 + c_2c_4+ c_1^2\,z) \nonumber\\
&&+2\,r^2\,A^{-1}\,\arctan\Big(\frac{\cos\theta-1}{\sin\theta}\Big)(c_2c_3-c_1c_4)
= d_3. \nonumber
\end{eqnarray}
when we assumed that
\begin{eqnarray*}
&&\hspace{-0.7cm} A = c_2^2\,z^2 + 2\,c_2c_4\,z + c_4^2 +
c_1^2\,z^2+2\,c_1c_3\,z +
c_3^2\\
&&\hspace{-0.7cm} B = (-2\,c_1c_3\,z-c_3^2 - c_1^2\,z^2 + c_5^2\,r^2 - c_2^2\,z^2 - 2\,c_2c_4\,z-c_4^2)^{\frac{1}{2}}\\
&&\hspace{-0.7cm} C = \frac{1}{\cos\theta+1}(c_1\,z\,\cos\theta +
c_3\,\cos\theta
-c_4\sin\theta- c_2\,z\,\sin\theta + c_5\,r)\\
&&\hspace{-0.7cm} \alpha = \frac{1}{B\sin\theta}(c_1 z - c_1
z\cos\theta + c_3 - c_3\cos\theta -c_5 r+ c_5 r\cos\theta + c_2
z\sin\theta + c_4\sin\theta).
\end{eqnarray*}
and $d_i$~s ($i=1,2,3$) are some constants.
\par
The functions $I_1,I_2,I_3$ form a complete set of functionally
independent invariants of one-parameter group generated by
(\ref{eq:9}) (see \cite{Ol}).\par
According to the theorem of section 4.3.3 of \cite{Ib}, the
derived invariants (\ref{eq:10}) as independent first integrals of
the characteristic system of the infinitesimal generators
(\ref{eq:7}), provide the general solution
$$S(r,\theta,z,u):=\mu(I_1(r,\theta,z,u),I_2(r,\theta,z,u),I_3(r,\theta,z,u)),$$
with an arbitrary function $\mu$.

%******************************************
\section{Classification of subalgebras}
It is well-known that the problem of classifying invariant
solutions is equivalent to the problem of classifying subgroups of
the full symmetry group under conjugation in which itself is
equivalent to determining all conjugate subalgebras \cite{Ol,Ov}.
The latter problem, tends to determine a list (that is called an
{\it optimal system}) of conjugacy inequivalent subalgebras with
the property that any other subalgebra is equivalent to a unique
member of the list under some element of the adjoint
representation i.e. $\overline{{\goth h}}\,{\rm Ad}(g)\,{\goth h}$
for some $g$ of a considered Lie group.

For finding one-dimensional subalgebras, this classification
problem is essentially the same as the problem of classifying the
orbits of the adjoint representation.\par
The adjoint action is given by the Lie series
\begin{eqnarray*}
{\rm Ad}(\exp(sX_i))X_j
=X_j-s[X_i,X_j]+\frac{s^2}{2}[X_i,[X_i,X_j]]-\cdots,
\end{eqnarray*}
where $s$ is a parameter and $i,j=1,\cdots,7$. Let $M_i(s):{\goth
g}\rightarrow{\goth g}$ defined by $X\mapsto {\rm
Ad}(\exp(s\,X_i))X$ be the adjoint representation on $\goth g$
(for $i=1,\cdots,7$ and small parameter $s$) and also $M_i(s)$ be
its matrix of representations with respect to the basis
$\{X_1,\cdots,X_7\}$, by an abuse of convention. After
straightforwardly computations, we find that
\begin{eqnarray*}
M_1(s)&=&  \left[
\begin{array}{cccccccc}
1&0&0&0&0&0&0\\0&1&0&0&0&0&0\\0&0&1&0&0&0&0\\0&0&0&\cos s&\sin
s&0&0\\0&0&0&-\sin s&\cos s&0&0\\0&0&0&0&0&\cos s&-\sin
s\\0&0&0&0&0&\sin s&\cos s
\end{array} \right],\\
M_2(s)&=& \left[
\begin{array}{cccccccc}1&0&0&0&0&0&0\\0&1&0&0&0&0&0\\0&0&1&0&0&0&0\\
0&0&0&1&0&0&0\\0&0&0&0&1&0&0\\0&0&0&0&-s&1&0\\0&0&0&-s&0&0&1
\end{array} \right],\\
M_3(s)&=&{\rm Id}_7,\\
M_4(s)&=&  \left[
\begin{array}{cccccccc} 1&0&0&0&-s&0&0\\0&1&0&0&0&0&0\\0&0&1&0&0&0&0\\
0&0&0&1&0&0&0\\0&0&0&0&1&0&0\\0&0&0&0&0&1&0\\0&s&0&0&0&0&1
\end{array} \right],\\
M_5(s)&=&  \left[
\begin{array}{cccccccc} 1&0&0&s&0&0&0\\0&1&0&0&0&0&0\\0&0&1&0&0&0&0\\
0&0&0&1&0&0&0\\0&0&0&0&1&0&0\\0&s&0&0&0&1&0\\0&0&0&0&0&0&1
\end{array}
\right],\\
M_6(s) &=&  \left[
\begin{array}{cccccccc} \cos s&0&0&0&0&0&\sin s_6\\0&\cos s&0&0&\sin s&0&0\\0&0&1&0&0&0&0\\
0&0&0&1&0&0&0\\0&-\sin s&0&0&\cos s&0&0\\0&0&0&0&0&1&0\\-\sin
s&0&0&0&0&0&\cos s
\end{array} \right],\\
M_7(s)&=& \left[
\begin{array}{cccccccc} \cos s&0&0&0&0&-\sin s&0\\0&\cos s&0&\sin s&0&0&0\\0&0&1&0&0&0&0\\
0&-\sin s&0&\cos s&0&0&0\\0&0&0&0&1&0&0\\ \sin s&0&0&0&0&\cos s&0\\
0&0&0&0&0&0&1
\end{array} \right].
\end{eqnarray*}
In the matrices $M_1,M_6$ and $M_7$, there are  minor matrices
which indicate 3-dimensional rotations, that are related to the
section ${\goth s}{\goth o}(3)$ of the Lie structure
(\ref{eq:8}).\par
\paragraph{Theorem 4.}{\em An optimal system of one-dimensional Lie subalgebras of (\ref{eq:1}) is provided by those
generated by\\[2mm]
\begin{tabular}{rlrl}
1)&$X_1$,              &\hspace{2cm}10)&$a\,X_4 + b\,X_6$,\\
2)&$X_2$,              & 11)&$a\,X_5 + b\,X_7$,\\
3)&$X_6$,              &12)&$X_3 + a\,X_1 + b\,X_2$,\\
4)&$X_7$,              &13)&$X_3 + a\,X_4 + b\,X_6$,\\
5)&$X_3 + a\,X_1$,     &14)&$a\,X_1 + b\,X_2 + c\,X_5$, \\
6)&$X_3 + a\,X_4$,     &15)&$a\,X_2 + b\,X_5 + c\,X_7$,\\
7)&$X_3 + a\,X_6$,     &16)&$a\,X_4 + b\,X_5 + c\,X_6$,\\
8)&$X_3+a\,X_7$,       &17)&$X_3+a\,X_1+b\,X_2 + c\,X_5$,\\
9)&$a\,X_1+b\,X_2$,
\end{tabular}\\[2mm]
where the coefficients $a$, $b$ and $c$ are arbitrary constants.}\\

\noindent{\it Proof:} Let $X=\sum_{i=1}^7 a_i\,X_i$ be an element
of ${\goth g}$, then we see that
\begin{eqnarray*}
&&\hspace{-1cm}M_7(s_7)\circ \cdots\circ M_2(s_2)\circ M_1(s_1)\;:\;X \;\longmapsto\; \Big(a_1\,\cos s_6 \cos s_7 +\cdots +a_7\,\sin s_6\Big)\,X_1\\
&+&\Big(a_2\,\cos s_6\cos s_7+a_4\,\cos s_6\sin s_7 +a_5\,\sin s_6\Big)\,X_2 \\[-1mm]
&\vdots&\\[1mm]
&+& \Big(a_1\,(\sin s_1\sin s_7-\cos s_1\sin s_6\cos
s_7)+\cdots+a_7\,\cos s_1\cos s_6\Big)\,X_7.
\end{eqnarray*}
For classifying one-dimensional Lie subalgebras of CHE, we plan
the following cases that in each case, by acting a finite number
of the adjoint representations $M_i$ ($i=1,\cdots,7$) on $X$, by
proper selection of parameters $s_i$ in each stage, we gradually
try to make the coefficients of $X$ vanish and to acquire the most
simple form of $X$.

$\bullet$ At first, we suppose that $a_3\neq 0$ and scaling if
necessary, we can assume that $a_3=1$. This assumption suggests
{\it Case 1,$\cdots$, Case 10}.

{\it Case 1.}  If $a_5\neq 0$, then we act on $X$ by ${\rm
Ad}(\exp(-\arctan(a_4/a_5)\,X_1))$ and hence we can make the
coefficient of $X_4$ vanish. Then we tend to the new form
$$X'=a_1\,X_1 + a_2\,X_2 + X_3 + a'_5\,X_5 + a'_6\,X_6 + a'_7\,X_7$$
for certain scalars $a'_5$ and $a'_6,a'_7$ depending on $a_4,a_5$
and $a_4$, $a_5$, $a_6$, $a_7$ resp. By acting ${\rm
Ad}(\exp(\frac{a'_7}{a_2}\,X_4))$ on $X'$ when we suppose that
$a_2\neq 0$, we tend to
$$X''=a'_1\,X_1 + a_2\,X_2 + X_3 + a'_5\,X_5 + a'_6\,X_6$$
where the coefficient of $X_5$ is vanished for scalar $a'_1$
depending to $a_1$, $a_2$, $a'_5$, $a'_7$. Also be the action of
${\rm Ad}(\exp(-(a'_6/a_2)\,X_5))$ on $X''$ the coefficient of
$X_6$ will be zero:
$$X'''=a'_1\,X_1 + a_2\,X_2 + X_3 + a'_5\,X_5.$$
At this stage, by acting adjoint representations $M_i(a_i)$ on
$X'''$, we find that no more simplification of $X'''$ is possible.
Thus each of $a'_1$ and $a'_5$ are arbitrary. Therefore, any
one-dimensional Lie subalgebra generated by $X$ with $a_2,a_5\neq
0$ (by knowing that $a_3\neq 0$) is equivalent to the Lie
subalgebra spanned by $$ a'_1\,X_1 + a_2\,X_2 + X_3 + a'_5\,X_5,$$
which is equal to section ${\it 17)}$ of the theorem.

\medskip Each of the following cases is prepared by a similar method to
{\it Case 1} and by eliminating unnecessary details, we just give
the conditions.

\medskip \noindent{\it Case 2.} If $a_5\neq 0$ and $a_2=0$, then we can
make the coefficients of $X_6,X_1$ vanish resp. by $M_2,M_4$ and
then $X$ is reduced to ${\it 11)}$.

\medskip \noindent{\it Case 3.} In the case which $a_5=0$ along with
$a_2,a_4\neq 0$, we act on $X$ by $M_4,M_5,M_5,M_7$ to cancel the
coefficients of $X_7,X_6,X_1,X_2$ resp., then the simplest
possible form of $X$ is equal to the part ${\it 6)}$ of the
theorem.

\medskip \noindent{\it Case 4.} When we change the condition $a_4\neq 0$ in
{\it Case 3.} to $a_4=0$, by applying $M_4,M_5$ the coefficients
of $X_7,X_6$ vanish and no more simplification is possible.
Therefore $X$ is reduced to ${\it 12)}$.

\medskip \noindent{\it Case 5.} When $a2=a_5=0$ and $a_4,a_7\neq 0$, then
by applying $M_4,M_6,M_2$ one can make the coefficients of $X_5$,
$X_1$, $X_7$ vanish and hence $X$ turns to ${\it 13)}$.

\medskip \noindent{\it Case 6.} If in {\it Case 5} we only change $a_4\neq
0$ to $a_4=0$ along with $a_6\neq 0$, then we can act on $X$ by
resp. $M_4$, $M_6$, $M_1$ to cancel the coefficients $X_5$, $X_1$,
$X_7$. Hence we find that form of $X$ is in the form of ${\it
7)}$.

\medskip \noindent{\it Case 7.} In the latter case, if we suppose that
$a_6=0$, by acting $M_4,M_6$ the coefficients of resp. $X_5$,
$X_1$ are zero and $X$ is equal to the part ${\it 8)}$ of the
theorem.

\medskip \noindent{\it Case 8.} Let $a2=a_5=a_7=0$ and $a_4\neq 0$. we can
make the coefficients of $X_5,X_1$ resp. by $M_4$, $M_5$ vanish
and then reduce $X$ to the section ${\it 13)}$.

\medskip \noindent{\it Case 9.} When $a2=a_4=a_5=a_7=0$ and $a_6\neq 0$, by
acting resp. $M_4,M_7$ on $X$, one can make the coefficients of
$X_5$, $X_1$ vanish. With these conditions, $X$ gives the form of
${\it 7)}$.

\medskip \noindent{\it Case 10.} When $a2=a_4=a_5=a_6=a_7=0$, by the action
of $M_4$ on $X$ and vanishing $X_5$, $X$ is reduced to ${\it 5)}$.

\medskip $\bullet$ At the second, we suppose that $a_3=0$ and consider {\it
Case 11,$\cdots$, Case 20} prepared by this assumption.

\medskip \noindent{\it Case 11.} Let $a_4,a_5\neq0$. We can act by $M_2$,
$M_7$, $M_4$ to cancel the coefficients of $X_7$, $X_2$, $X_1$
resp. then we will have the reduced form ${\it 16)}$ of $X$.

\medskip \noindent{\it Case 12.} By changing $a_5\neq0$ of the latter case
to be zero and applying $M_2$, $M_7$, $M_5$ on $X$, the
coefficients of $X_7$, $X_2$, $X_1$ resp. will vanish and then $X$
is reduced to ${\it 10)}$.

\medskip \noindent{\it Case 13.} Let $a_4=0$, $a_7\neq0$ and $a_5\neq0$ (by
following the process and finding the exact generating
coefficients, in fact, we assume that $a'_5 = a_1 a_2+ a_5 a_7$ be
non-zero). The coefficients of $X_1$, $X_6$ by acting $M_6$, $M_2$
will be zero and no further simplification is possible. Hence $X$
is reduced to ${\it 15)}$.

\medskip \noindent{\it Case 14.} Suppose that $a_4=a'_5=0$, $a_7\neq0$, and
$a_2\neq 0$ (in fact, in the method of construction, $a'_2=-a_2
a_7 + a_1 a_5\neq0$). Applying $M_6,M_4,M_5$ resp. on $X$ make the
coefficients of $X_1$, $X_7$, $X_6$ vanish. In this case $X$ is
reduced to ${\it 2)}$.

\medskip \noindent{\it Case 15.} In the last case, we let $a'_2$ be zero
and apply $M_6$, $M_1$, $M_1$ resp. to cancel the coefficients of
$X_1$, $X_2$, $X_6$. Then $X$ will be in the form of ${\it 4)}$.

\medskip \noindent{\it Case 16.} If $a'_2=a_4=a'_5=b_7=0$ then by acting
$M_6$, $M_1$ resp. to cancel the coefficients of $X_1$, $X_6$, we
lead to the simplest case ${\it 3)}$ of $X$.

\medskip \noindent{\it Case 17.} Let $a_4=a_7=0$ and $a_5\neq 0$. By acting
$M_2$ on $X$ then the coefficient of $X_6$ will be zero and hence
$X$ is reduced to ${\it 14)}$.

\medskip \noindent{\it Case 18.} If $a_4=a_5=a_7=0$ and $a_2\neq 0$, then
by acting $M_5$ to cancel the coefficient of $X_6$, we lead to the
simplest case ${\it 9)}$ of $X$.

\medskip \noindent{\it Case 19.} When $a_2=a_4=a_5=a_7=0$ and $a_6\neq 0$,
then we can make the coefficients $X_1$ resp. by $M_7$ vanish.
Then $X$ is reduced to ${\it 3)}$.

\medskip \noindent{\it Case 18.} Finally, if in the last case $a_6=0$, no
further simplification is possible and then $X$ is reduced to
${\it 1)}$.

\medskip \noindent There is not any more possible case for studying and the
proof is complete.\hfill\ $\square$

\medskip For finding optimal systems of higher dimensional Lie
subalgebras (equivalently Lie subgroups), we refer the reader to
\cite{Ov} or \cite{Ib1}.

\end{document}